\documentclass
[12pt,a4paper]{article}
\usepackage[cp1251]{inputenc}
\usepackage{amssymb, amsfonts, amsmath, latexsym}
\usepackage[english]{babel}

\begin{document}


\Large\centerline{\bf  Selective survey on spaces of closed subgroups of topological groups  }\vspace{6 mm}

\normalsize\centerline{\bf  Igor V. Protasov}\vspace{6 mm}

{\bf Abstract.}  We survey   different topologizations of the set $\mathcal{S}(G)$ of all closed subgroups of a topological group $G$  and demonstrate some applications in {\it Topological  Grous, Model Theory, Geometric Group Theory, Topological Dynamics.}
\vspace{6 mm}


\vspace{3 mm}

\vspace{6 mm}



Some words in place of introduction. For a topological group $G$,  $\mathcal{S}(G)$  denotes the set of all closed subgroup of $G$.
There are many ways to endow $\mathcal{S}(G)$  with a topology  related  to the topology of $G$.
Among them, the most intensively studied    are Chabauty topology rooted in {\it Geometry of Numbers}  and the Vietories topology went from {\it General Topology};  both coincide if $G$  is  compact.
The spaces of closed subgroups are interesting by their own  sake but also have   some deep applications in {\it Topological  Groups and Model Theory, Geometric Group Theory and Dynamical Systems.}
The survey is my subjective look at this area.

\vspace{3 mm}
Content: Chabauty spaces;  Victories spaces;  Other topologizations.
\vspace{3 mm}

\section{Chabauty spaces}

{\bf 1.1.   From Minkowski  to Chabauty.}
We recall that a {\it lattice} $L$ in  $\mathbb{R}^{n}$  is a discrete subgroup of rank $n$.
 We denote $min \  L$ the length of a shortest  non-zero vector from $L$, $vol \ (\mathbb{R}^{n}/ L)$  is the volume of a basic parallelepiped of $L$ .

A sequence $(L_{m})_{m \in\omega} $ of lattices in $\mathbb{R}^{n}$ converges to the lattice $L$ if,  for each $m\in \omega$,  one can choose a basis
$a_{1} (m), \ldots,  a_{n} (m)$  of $L_{m}$ and a basis
 $a_{1}, \ldots ,   a_{n}$  of $L$ such  that the sequence $(a_{i}(m)) _{m\in\omega} $ converges to $a_{i}$  for each $i\in \{1, \ldots, n\}$.  This convergence of lattices was introduced by  H. Minkowski  \cite{b1}, and its  usage in  {\it Geometry of  Numbers} (see \cite{b2}) is based on the following theorem of K.  Mahler \cite{b3}.
\vspace{5 mm}

{\bf  Theorem 1.1.} {\it  Let $\mathcal{M}$ be a set of lattices in  $\mathbb{R}^{n}$.
Every  sequence in $\mathcal{M}$  has a convergent subsequence if and only if  there exist  two constants $C>0$,  $c>0$  such that $min\  L > c$, $vol \ (\mathbb{R}^{n}\setminus   L) < C$   for each  $L\in \mathcal{M}$}.
\vspace{5 mm}

What we know now as Chabauty topology was invented by C. Chabauty in \cite{b4} in order to extend Theorem 1.1 to lattices in connected Lie groups.
A discrete subgroup $L$ of a connected Lie group $G$   is called a {\it lattice} if the quotient space $G/  L$ is compact.

Let $X$ be a Hausdorff locally compact space and let $exp \ X$ denotes the set  of all closed subsets of $X$.
The sets
$$\{ F\in exp \  X: F\cap K = \emptyset\},  \  \{F\in exp \  X: F\cap U \neq \emptyset\},$$
where $K$ runs over all compact subsets of $X$  and  $U$ runs over all open subsets of $X$ ,  form the subbase of the {\it Chabauty topology} on $exp \  X$.
The   space $exp \  X$ is compact and Hausdorff.  If $X$  is discrete then $exp \  X$ is homeomorphic to the Cantor cube
$\{ 0, 1\} ^{|X|}$.

We note also that a net $( F_{\alpha} )_{\alpha\in\mathcal{I}}$  converges in $exp \  X$ to $F$  if and only if

\begin{itemize}
\item{} for every compact $K$ of $X$ such that $K\cap F = \emptyset$,  there exists
$\beta\in\mathcal{I}$   such that  $F_{\alpha}\cap K=\emptyset$
 for each $\alpha>\beta$; \vskip 5pt

\item{}  for every  $x\in F$ and   every  neighbourhood $U$  of $x$,   there exists
$\gamma\in \mathcal{I}$  such that $F _{\alpha} \cap U\neq\emptyset$  for each $\alpha>\gamma$. \end{itemize}

If $G$ is a locally compact group then $\mathcal{S}(G)$ is a closed subspace of $exp \  G$  (so  $\mathcal{S}(G)$ is compact);    $\mathcal{S}(G)$ is called the {\it Chabauty space} of $G$.
\vskip 6pt

{\bf Theorem 1.2[4]. }
{\it Let $G$ be a connected unimodular Lie  group.
A set $\mathcal{M}$ of lattices in $G$ is relatively compact in $\mathcal{M}$
 if and only if there exists constant $C> 0$  and a neighbourhood $U$  of the identity $e$ of $G$ such that $L\cap U = \{e\}$ and $vol  \ (G/ L)< C$ for each  $L\in  \mathcal{M}$.}
\vskip 5pt

With some technical improvement made in [5], the paper [4]is included in [6, Chapter 8].
\vskip 7pt

{\bf 1.2. Pontryagin-Chabauty duality.}
This duality was established in \cite{b7}  and detalized in \cite{b8}.
We use the standard abbreviation LCA for a locally compact Abelian group.  Let $G$ be a LCA-group, $G^{\wedge}$
 denotes its dual group, $G^{\wedge}= Hom \ (G, \mathbb{R} /\mathbb{Z})$
   and let $\varphi$ denotes the canonical bijection $\mathcal{S}(G)\longrightarrow \mathcal{S}(G^{\wedge})$, $\varphi(X)= \{f\in G^{\wedge}: X\subseteq ker \  f\}$.
   \vskip 5pt

{\bf  Theorem 1.3.}  {\it For every  LCA-group  $G$, the bijection $\varphi: \mathcal{S}(G)\longrightarrow \mathcal{S}(G^{\wedge})$ is a homeomorphism.} \vskip 5pt

Typically, Theorem 1.3 applies to replace $\mathcal{S}(G)$
 by $\mathcal{S}(G^{\wedge})$
  in the case of a compact Abelian group $G$.

In what follows we use the notations:  $\mathbb{C}_{n}$  is the cyclic group of order   $n$ ,
$\mathbb{C}_{p^{\infty}}$  is the quasi-cyclic (or Pr$\ddot{u}$ffer)
$p$-group,  $\mathbb{Z}$ is the discrete group of integers,  $\mathbb{Z}_{p}$  is the group of $p$-adic  integers,
 $\mathbb{Q}_{p}$  is the additive group of the field of  $p$-adic   numbers.

\vskip 7pt

{\bf 1.3.   $S(G)$ for compact $G$.}
The following two lemmas from \cite{b9}   are the basic technical  tools in this area.
\vskip 5pt

{\bf Lemma  1.1.}  {\it If $G,H$ are compact groups and $\varphi: G \longrightarrow H$ is a continuous surjective  homomorphism then the mapping $\mathcal{S}(G)\longrightarrow  \mathcal{S}(H),$   $X\longmapsto  \varphi (X)$ is continuous and open.}\vskip 5pt

The continuity is easy but to  prove  the openness we need
\vskip 5pt

{\bf Lemma  1.2.}  {\it  Let $G$  be a compact group, $X\in \mathcal{S}(G)$.
Then the following subsets from a base of neighbourhoods  of $X$ is $\mathcal{S}(G)$:

$$\mathcal{N}_{X}(U, N,x_{1}, \ldots , x_{n})=  \{u^{-1}Y u: Y\in\mathcal{S}(G),  \  Y\subseteq XN, \  Y\cap x_{1} U \neq \emptyset , \ldots , Y \cap x_{n}U\neq\emptyset , \}$$
 where  $U$  is a neighbourhood of the identity of $G$,  $N$  is closed normal subgroup such that $G/N$ is a Lie group, $x_{1},\ldots,  x_{n}$  are arbitrary  elements of} $X$,  $n\in \mathbb{N}$.
\vskip 5pt

In particular,  if $G$  is a compact Lie group then  Lemma  1.2  states that there is a  neighbourhood  $\mathcal{N}$ of $X$  such that each subgroup $Y\in \mathcal{N}$   is conjugated to some subgroup of $X$.
The key part in the proof of   Lemma  1.2 plays the Montgomery-Yang theorem on tubes \cite{b10}, see also [11, Theorem 5.4 from Chapter 2].

We recall that the {\it cellularity}  (or Souslin number) $c(X)$  of a topological space $X$  is the supremum  of cardinalities of disjoint families of open subsets of $X$. A topological space $X$  is called {\it dyadic}  if $X$  is a continuous image of some Cantor cube $\{0,1\}^{\kappa}$.

The {\it weight}   $w(X)$   of a topological space $X$  is the minimal cardinality of open bases of $X$.
\vskip 7pt

{\bf Theorem 1.4 [9].  } {\it For every compact group $G$, we have $c(\mathcal{S}(G))\leq \aleph_{0} $.
In addition, if $w(G)\leq \aleph_{1} $  then  $\mathcal{S}(G)$ is dyadic. }
\vskip 7pt

{\bf Theorem 1.5 [12].  } {\it
Let a group $G$ be either profinite or compact and Abelian.  If $w(G)> \aleph_{2}$  then the space  $\mathcal{S}(G)$  is not dyadic.}
\vskip 7pt

{\bf Theorem 1.6 [12].  } {\it Let $G$  be an  infinite compact Abelian group such that $w(G)\leq\aleph_{1}$.
Then the  space  $\mathcal{S}(G)$ is homeomorphic to the Cantor cube $\{0,1\}^{w(G)}$  if and only if $\mathcal{S}(G)$  has no isolated points.}
\vskip 5pt

An Abelian group $G$ is called {\it Artinian} if every increasing chain of subgroups of $G$  is finite; every such group is isomorphic to the direct $sum \ \oplus_{p\in F} \mathbb{C}_{p^{\infty}}\oplus K$,  where $F$  is a finite set  of primes,  $K$ is a finite subgroup. An Abelian  group $G$  is called {\it minimax}  if $G$ has a finitely generated subgroup $N$  such that $G / N$ is Artinian.
\vskip 7pt

{\bf Theorem 1.7 [12].  } {\it  For a compact Abelian group $G$,  the space  $\mathcal{S}(G)$ has an isolated point if and only if the dual group  $G^{\wedge}$  is minimax.}
\vskip 7pt

{\bf 1.4. $S(G) $  for LCA $G$.  }  The space $\mathcal{S}(\mathbb{R})$  is homeomorphic to the segment $[0,1]$.
By [13],  $\mathcal{S}(\mathbb{R} ^{2})$ is homeomorphic to the sphere  ${\bf S }^{4}$.   For $n\geq  3$, $\mathcal{S}(\mathbb{R} ^{n})$  is not a topological manifold and its  structure is   far from understanding, see [14].
\vskip 7pt

{\bf Theorem 1.8 [15].  } {\it  The  space $\mathcal{S}(G) $   of a LCA-group $G$ is connected if and only if $G$ has a subgroup topologically isomorphic to $\mathbb{R}$.}

If $F$ is a non-solvable  finite group then $\mathcal{S}(\mathbb{R} \times  F)$  is not connected [8, Proposition 8.6].

\vskip 7pt

{\bf Theorem 1.9 [8].  } {\it  The  space $\mathcal{S}(G) $   of a LCA-group $G$ is
totally disconnected if and only if $G$
is either totally disconnected or each elements of $G$
belongs to a compact subgroup.}
\vskip 7pt

Some more  information on $\mathcal{S}(G) $   for LCA $G$  can be find in [8] and references there, in particular, on topological dimension of
$\mathcal{S}(G) $.

By Theorems 1.4 and 1.3,  $c(\mathcal{S}(G))\leq\aleph_{0} $ for every discrete Abelian group.  We say that a  topological space $X$  has {\it Shanin number} $\omega$   if any  uncountable  family $\mathcal{F}$ of non-empty open subsets of $X$ has an uncountable subfamily   $\mathcal{F}^{\prime}$  such that $\cap\mathcal{F}^{\prime}\neq\emptyset$ .
Evidently, if a space $X$ has Shanin number $\omega$ then  $c(X)\leq\aleph_{0}$ .
By [16, Theorem 1], for every discrete Abelian group $G$, the space
 $\mathcal{S}(G) $ has Shanin number $\omega$.
 By [16, Theorem 3], for every infinite cardinal $\tau$, there exists a solvable discrete group $G$ such that $c(\mathcal{S}(G))=|G|=\tau$.

\vspace{7 mm}

{\bf 1.5. $S(G) $ as a lattice.  }  The  set $S(G) $  has the natural structure of a lattice with the operations
$\vee$  and $\wedge$,
 where $A \wedge B= A\cap B$  and $A\vee B$  is the smallest closed subgroup of $G$ containing $A$ and $B$.
 In this subsection, we formulate    some results from [17] on interrelations between the topological and lattice structures on $S(G) $.

For $g\in G$,   $\overline{<g>}$ denotes the  subgroup of $G$ topologically generated by $g$.
A totally disconnected locally compact group $G$ is called {\it periodic}  if $\overline{<g>}$
 is compact for each  $g\in G$
In this case, $\pi(G)$  denotes the set of all prime  numbers such that $p\in\pi(G)$
  if and only if there is  $g\in G$ such that $\overline{<g>}$
    is topologically isomorphic either to $\mathbb{C}_{p^{n}}$
     or to $\mathbb{Z}_{p}$;  this $g$ is called a {\it topological $p$-element}.
\vspace{5 mm}

{\bf Theorem 1.10.  } {\it  For a compact group $G$,  the following statements are equivalent
\vspace{3 mm}

$(i)$  $\wedge$ is continuous;
\vspace{3 mm}

$(ii)$ $\wedge$ and $\vee$ are  continuous;
\vspace{3 mm}

$(iii)$  $G$ is the semidirect product $K \leftthreetimes P$,  where $K$ is profinite  with finite Sylow $p$-subgroups, $P$ is Abelian profinite and each Sylow $p$-subgroup of $G$ is $\mathbb{Z}_{p}$,   $\pi(K)\cap \pi(P)=\emptyset$ and the centralizer of each Sylow $p$-subgroup of $G$  has finite index in $G$.}
\vspace{5 mm}

{\bf Theorem 1.11.  } {\it  For a locally compact group $G$,  the operation $\wedge$ is continuous if and only if the
 followings  conditions are satisfied

 \vspace{3 mm}

$(i)$  $G$ is either discrete or periodic;
\vspace{3 mm}

$(ii)$ $\wedge$ is continuous in $\mathcal{S}(H)$ for  each compact subgroup $H$ of $G$;
\vspace{3 mm}

$(iii)$  the centralizer of each topological $p$ element of $G$ is open.
}
\vspace{5 mm}

We recall that a torsion group $G$ is {\it layerly  finite} if the set
 $\{g\in G: g^{n}=e\}$  is finite for each $n\in \mathbb{N}$.
A layerly finite group $G$ is called {\it thin}  if each Sylow $p$-subgroup of $G$ is finite (equivalently, $G$ has no subgroup isomorphic to $\mathbb{C}_{p^{\infty}}$).
\vspace{7 mm}

{\bf Theorem 1.12.  } {\it Let $G$ be a locally compact group. The operations
$\wedge$  and $\vee$ are continuous if and only if $G$ is periodic and topologically isomorphic  to
 $A\times B\times (C\leftthreetimes D)$,  where
  $C$  has a dense thin layerly  finite  subgroup, $A, B, D$ are Abelian  with Sylow $p$-subgroups $\mathbb{C_{p^{\infty}}}$, $\mathbb{Q}_{p}$  or $\mathbb{Z}_{p}$,
   the sets  $\pi(A)$, $\pi(B)$, $\pi(G)$, $\pi(D)$   are  pairwise disjoint and the centralizer of each Sylow  $p$-subgroup of $G$ is open.}
\vspace{7 mm}

{\bf 1.6. From Chabauty to local method.} A topological group $G$ is called {\it topologically simple}  if each closed  normal subgroup  of  $G$  is either $G$  or $\{e\}$.
Every  topologically  simple LCA-group is discrete and either  $G= \{e\}$  or  $G$ is isomorphic to $\mathbb{C}_{p}$.

Following the algebraic tradition,  we say that a group  $G$ is {\it locally nilpotent (solvable)}  if every finitely generated subgroup is nilpotent (solvable).

In [18, Problem 1.76],  V. Platonov  asked whether there exists  a non-Abelian topologically simple locally compact locally nilpotent group. Now we sketch the negative answer to this  question for locally solvable group obtained in [19].

Let $G$  be a locally compact locally solvable group.  We   take  $g\in G \ \setminus \{e\}$,  choose a compact neighbourhood  $U$ of $G$ and denote by $\mathcal{F}$ the family of  all topologically finitely generated subgroups of $G$ containing $g$.
We may assume that $G$  is not topologically finitely
generated so $\mathcal{F}$  is  directed by the inclusion $\subset$.
For each  $F\in \mathcal{F}$, we choose  $A_{F}$, $B_{F} \in \mathcal{S}(F)$  such that
$B_{F}\subset  A_{F}$,  $A_{F}$  and $B_{F}$ are normal in $F$, $A_{F}\cap U\neq\emptyset$,  $B_{F}\cap U=\emptyset$ and  $A_{F} /  B_{F}$  is Abelian. Since $\mathcal{S}(G)$  is compact,  we can choose two  subnets
$(A_{\alpha})_{\alpha\in \mathcal{I}}$, $(B_{\alpha})_{\alpha\in \mathcal{I}}$  of the nets
$(A_{F})_{F\in \mathcal{F}}$, $(B_{F})_{F\in \mathcal{I}}$
 which  converges to $A, B\in \mathcal{S}(G)$.
Then $A, B$  are normal in $G$  and   $A/B$  is Abelian. Moreover, $x \notin B$ and  $A\cap U\neq\emptyset$.
If $A \neq \{G\}$  then $A$ is a proper normal subgroup of $G$;
 otherwise $G/B$ is Abelian.

\vspace{3 mm}

In [20], the Chabauty     topology was  defined on some systems of closed  subgroups  of                                              locally compact group $G$.
A system $\mathfrak{A}$ of  closed subgroups  of  $G$ is  called  {\it subnormal}  if

\begin{itemize}

\item{} $\mathfrak{A}$ contains $\{e\}$ and $G$;

\item{} $\mathfrak{A}$ is linearly ordered by the inclusion $\subset$;

\item{} for any subset $\mathfrak{M}$ of $\mathfrak{A}$,  the closure of
$\bigcup_{F\in\mathfrak{M}}\ F\in \mathfrak{A}$   and $\bigcap_{F\in\mathfrak{M}}\ F\in \mathfrak{A}$ ;

\item{} whenever $A$  and $B$ comprise  a jump in $\mathfrak{A}$ (i.e $B\subset A$  and  no members of $\mathfrak{A}$ lie between  $B$  and $A$),  $B$  is a normal subgroup of $A$.

\end{itemize}

If the subgroup $A, B$ form a jump then  $A/ B$ is called  a factor of  $G$.
The system is called {\it normal}  if  each $A\in \mathfrak{A}$  is normal in $G$.

A group   $G$ is called an RN-group if $G$ has a normal
system with Abelian factors. Among the local theorems
  from  [20],  one can find the following: if every topologically
   finitely    generated  subgroup of a locally compact group  $G$ is an  RN-group then $G$  is an RN-group.
In particular, every locally compact  locally solvable  group is an RN-group.

In 1941, see [21, pp. 78-83], A.I. Mal'tsev obtained  local  theorems
 for discrete groups as applications of the following general local theorem: if every   finitely    generated  subsystem of an algebraic system $A$  satisfies some property $\mathcal{P}$,  which can be defined by some quasi universal second order  formula,  then $A$ satisfies  $\mathcal{P}$.

In [22], Mal'tsev's  local theorem was  generalized on topological algebraic system.
The part of the
 model-theoretical  Compactness Theorem in  Mal'tsev arguments  plays  some  convergents of  closed subsets.
  A net  $(F_{\alpha})_{\alpha\in\mathcal{I}}$ of  closed  subsets  of a topological space $X$    $S$-converges to  a closed subset  $F$ if

\begin{itemize}

\item{}  for every  $x\in F$ and every neighbourhood $U$  of $x$, there exists  $\beta\in \mathcal{I}$
 such  that $F_{\alpha}\cap U\neq\emptyset$  for each $\alpha>\beta$;

\item{}  for every $y\in X\setminus F$,  there exist a neighbourhood $\mathcal{V}$  of $y$  and   $\gamma\in\mathcal{I}$ such that $F_{\alpha}\cap \mathcal{V}= \emptyset$   for each  $\alpha>\gamma$.
\end{itemize}

Every net of closed subsets of an arbitrary (!)
 topological space has a convergent subnet.
If $X$ is a Hausdorff  locally compact space then $S$-convergence  coincides with  convergence in the Chabauty topology.
\vspace{7 mm}

{\bf  1.7  Spaces of marked groups}.  Let $F_{k}$ be the free group of rank $k$  with  the  free generators $x_{1},\ldots , x_{k}$   and let $\mathcal{G} _{k}$  denotes the set of all normal subgroups of $F_{k}$.   In the metric form, the Chabauty  topology on  $\mathcal{G}_{k}$   was introduced in [23] as a reply on the Gromov's idea of topologizations   of some sets of groups [24].

Let $G$ be a group generated by $g_{1},\ldots , g_{k}$. The bejection $x_{i}\longmapsto g_{i}$  $g_{1},\ldots , g_{n}$
 can be extended to the homomorphism $f: F_{k}\longrightarrow   G$.
With the correspondence   $G\longmapsto  ker \  f$, $\mathcal{G}_{k}$  is called the  {\it space marked $k$-generated groups}.

A couple of papers in development of  [23]  is directed to understand  how large in topological sense are well-known classes of finitely generated groups, or how a given  $k$-generated group is  placed in $\mathcal{G}_{k}$, see [25].
Among applications of $\mathcal{G}_{k}$,   we mention the construction of topologizable Tarski Monsters in [26].
\vspace{7 mm}

{\bf  1.8  Dynamical development.}
Every locally compact group $G$  acts on the Chabauty space $\mathcal{S}(G)$  by the rule:
 $(g,H)\longmapsto g ^{-1} H g$.
Under this action, every   minimal closed invariant subset of  $\mathcal{S}(G)$   is called a {\it uniformly  recurrent subgroup},  URS for short.
The study of URSs was initiated  by Glasner  and Weiss  [27] with the following observation.

Let a locally compact group $G$  acts on a compact $X$ so that is $G$ minimal,
 i.e. the orbit of each point  $x\in X$ is dense.  We consider the mapping $Stab:  X \longrightarrow \mathcal{S}(G) $  defined  by  $Stab  (x)=\{g\in G:  gx =x\}$.
Then there is the unique URS  contained in the closure of   $Stab  (X)$. This URS is called the {\it stabilizer URS}.
Glasner and Weiss asked whether every URS of a locally compact group $G$  arises as the stabilizer
URS of a minimal  action of $G$ on a compact space.
This question was answered   in the affirmative in [28].

\section{Vietoris spaces}

For a topological space $X$, the Victoris topology on the set $exp \ X$  of all closed subsets of $X$ is defined by the subbase  of open sets
$$\{F\in exp \ X: F\subseteq U\},  \  \{F\in exp \ X: F\cap       V\neq\emptyset\},$$
     where $U, \mathcal{V}$  run over all  open subsets of $X$.

A net $(F_{\alpha})  _{\alpha\in \mathcal{I}}$ converges to $F$  in
$exp \ X$   if and only if

\begin{itemize}

\item{}  for each open subset   $U$ of $X$  such   that $F\subseteq U$,  there exists $\beta\in\mathcal{I}$  such that   $F_{\alpha}\subseteq U$    for each $\alpha>\beta$;

\item{}   for each $x\in F$ and each neighbourhood $\mathcal{V}$  of  $x$,
 there exists $\gamma\in\mathcal{I}$  such that
 $F_{\alpha}\cap \mathcal{V}\neq\emptyset$  for each $\alpha>\gamma$.
\end{itemize}

\vspace{3 mm}

If $X$  is regular then $\mathcal{S}(G)$  is  closed in $exp \  G$.  As
   to my knowledge, the spaces $\mathcal{S}(G)$,
     where $G$  needs  not  to be compact, endowed with the Vietoris  topologies
     appeared  in  [29] with characterization  of  LCA-groups  $G$
     such that the canonical  mapping $\varphi: \mathcal{S}(G)\longrightarrow \mathcal{S}(G^{\wedge})$
      is a homeomorphism.

\vspace{7 mm}
{\bf 2.1. Compactness. }
It is naively   to ask a constructive
description  of arbitrary topological groups $G$  with
 compact   space $\mathcal{S}(G)$  because  we know  nothing  even about  $G  $ with  $S(G)=2$.

\vspace{5 mm}
{\bf Theorem 2.1.}   [30].
{\it  Let $G$  be a locally compact group.  The space $\mathcal{S}(G)$   is compact if and only if $G$  is one of the following groups
\vspace{3 mm}

  $(i)$   $G$  is compact;
  \vspace{3 mm}

  $(2)$  $\mathbb{C}_{p_{1}^{\infty}}\times \ldots\times \mathbb{C}_{p_{n}^{\infty}}\times K$,   where  $p_{1},\ldots ,  p_{n}$ are distinct prime numbers, $K$  is finite and each $p_{i}$  is not a divisor of $|K|$;
\vspace{3 mm}

$(3)$  $Q_{p}\times K$,  where  $K$ is finite and  $p$  does not divide  $|K|$.}
\vspace{5 mm}

Similar characterization of groups with compact  $\mathcal{S}(G)$ is given in [31] provided  that $G$  has a base at the identity consisting of subgroups.
\vspace{7 mm}

{\bf Theorem 2.2.}  [32]. {\it Let $G$ be a locally compact group. A  closed subset  $\mathcal{F}$ of  $\mathcal{S}(G)$  is compact if and only  if the following conditions are satisfied
\vspace{5 mm}

$(i)$   every descending chain of  non-compact   subgroups  from $F$ is finite;
\vspace{5 mm}

$(ii)$   every closed subset  $\mathcal{F}^{\prime}$  of $\mathcal{F}$ has only finite number of non-compact subgroups maximal in $\mathcal{F}$;
\vspace{5 mm}

$(iii)$   if a closed subset  $\mathcal{F}^{\prime}$  of $\mathcal{F}$ has no non-compact subgroups then $\cup \mathcal{F}^{\prime}$ is compact.}

\vspace{7 mm}

Two corollaries:  every compact in $\mathcal{L}(G)$ consisting of  non-compact subgroups is scattered; a  subset  $\mathcal{F}$ is compact if and only if $\mathcal{F}$  is countably compact.

For locally compact groups with $\sigma$-compact  space  $\mathcal{S}(G)$  see [33],  a description of LCA-groups with   locally compact space $\mathcal{S}(G)$  is obtained in [34].

A topological group  $G$  is called {\it inductively compact} if every finite subset of $G$  is contained in compact subgroup.
For a group $G$,   $K(G)$  and  $IK(G)$ denote the sets of all compact and closed  inductively compact subgroups.
\vspace{7 mm}

  {\bf Theorem 2.3.}  [35]. {\it For every locally compact group $G$,  $IK(G)$ is the closure of $K(G)$.}
\vspace{7 mm}

Two corollaries: if $G$  is a connected Lie group then $K(G)$  is closed;  $\mathcal{S}(G)$   is a $k$-space  for each locally compact group $G$ of countable weight,  i.e.  the topology  of  $\mathcal{S}(G) $ is  uniquely  determined by the family of all compact subsets   of $\mathcal{S}(G) $.
\vspace{7 mm}

{\bf 2.2. Metrizability and normality. }  LCA-groups  $G $  with metrizable and normal  space $\mathcal{S}(G)$
  are characterized by  S. Panasyuk in the candidate thesis  {\it Normality and metrizability  of the space of closed subgroups},  Kyiv University, 1989.
These lists are completely constructive but too cumbrous so we formulate only
\vspace{5 mm}

  {\bf Theorem 2.4.}   {\it For a discrete Abelian  group  $G$,  the following statements are equivalent
\vspace{3 mm}

$(i)$  $\mathcal{S}(G)$ is  metrizable;
\vspace{3 mm}

$(ii)$  $\mathcal{S}(G)$ is  normal;
\vspace{3 mm}

$(iii)$  $G$  has a finitely generated  subgroup $H$ such that
$G/ H=  \mathbb{C}_{p_{1}^{\infty}}\times \ldots\times \mathbb{C}_{p_{n}^{\infty}}$,   where  $p_{1},\ldots ,  p_{n}$ are distinct primes.}

\vspace{3 mm}

In  general case,  metrizability and normality of $\mathcal{S}(G)$ are not equivalent but if $G$ a
  connected semisimple  Lie group then $\mathcal{S}(G)$  is   metrizable iff
  $\mathcal{S}(G)$  is normal iff $G$  is compact,  see [36],  [37].
The space $\mathcal{S}(G)$  of every connected solvable   Lie group is  metrizable  [36].
\vspace{5 mm}

{\bf 2.3.  Some  cardinal  invariants}.
We remind that $c(X)$  denotes  the  cellularity of $X$.
\vspace{4 mm}

{\bf Theorem 2.5.[9]}   {\it For  every infinite locally compact group $G$, we have}  $c(\mathcal{S}(G))\leq c(G)$.
\vspace{4 mm}

{\bf Theorem 2.6. [38].}   {\it For every locally compact group $G$, the  following  conditions  are equivalent
\vspace{3 mm}

$(i)$  $\mathcal{S}(G)$ is of countable pseudocharacter; \vspace{3 mm}

$(ii)$  $\mathcal{S}(G)$ is of countable tightness;\vspace{3 mm}

$(iii)$  $\mathcal{S}(G)$ is sequential;}\vspace{3 mm}

$(iv)$  $w(G)\leq\aleph_{0}$.

\section{Other topologizations}

{\bf 3.1.  Bourbaki uniformities}.
Let $(X, \mathcal{U})$  be a uniform space.
The uniformity $\mathcal{U}$  induces  the uniformity
 $\widetilde{\mathcal{U}}$ on the set  $\mathcal{F}(X)$ all  non-empty  closed subsets  of $X$
  which  has as a base the  family  of  sets of  the form
$$\{(A,B)\in\mathcal{F}(X) \times\mathcal{F}(X): B\subseteq U(A), \  A\subseteq U(B)\},$$
whenever  $U\in \mathcal{U}$.
The uniformity $\widetilde{\mathcal{U}}$   was  introduced in [39, Chapter 2, $\S$ 1]   and $\widetilde{\mathcal{U}}$   is called {\it the Bourbaki}  (sometimes, Hausdorff-Bourbaki)  {\it uniformity}.

Let $G$ be a topological group.  We endow $G$  with the left uniformity $L$
 and  $F(G)$  with the  Bourbaki uniformity $\widetilde{L}$.
We denote by $\mathcal{L}(G)$  and  $\mathcal{B}(G)  $  the
subspaces of $\mathcal{F}(G)$  consisting of all   subgroups  and all totally  bounded subsets of $G$.
\vspace{5 mm}

{\bf Theorem 3.1.[40]}   {\it  Let a group $G$ has a base at the
identity consisting  of subgroups.  The space
$\mathcal{L}(G)$ is compact if  and only if  $G$  is totally bounded and
$K\bigcap  G$
 is dense in $K$  for each  closed  subgroup  $K$ from the completion of $G$.}

\vspace{5 mm}
 In particular, if   $\mathcal{L}(G)$  is compact then  $G$ is totally  minimal.
\vspace{5 mm}

{\bf Theorem 3.2.[40]}   {\it If a group $G$  is complete in the  left uniformity then
$\mathcal{B}(G)$  is  complete.}
\vspace{5 mm}

We recall that a topological group $G$  is {\it almost metrizable} if each neighbourhood of $e$ contains a
compact subgroup $K$ such that the set of all open  subsets  containing $K$  has a countable base.
Every metrizable  and every locally compact topological group are   almost metrizable.
\vspace{5 mm}

{\bf Theorem 3.3.[40]}   {\it
If an almost  metrizable  group  $G$  is complete in the left  uniformity then  $\mathcal{F}(G)$ is complete.}
\vspace{5 mm}

In [41],  Theorem 3.3  is proved with the bilateral  uniformity on $G$  (and so on $\mathcal{F}(G)$) in place of  the left uniformity.
\vspace{7 mm}

{\bf 3.2.
Functionally  balanced  groups.}  For a topological group  $G$,  the set $\mathcal{F}(G)$ has  the  natural  structure of a  semigroup with the operation  $(A, B)\longmapsto   cl \  AB$.
\vspace{5 mm}

{\bf Theorem 3.4.[42]}   {\it
For a topological  group $G$, the following  statements are  equivalent
\vspace{3 mm}

$(i)$  $\mathcal{F}(G)$ is a topological semigroup;

\vspace{3 mm}
$(ii)$  for every subset $X$ of $G$ and every neighbourhood $U$ of
$e$,  there exists a neighbourhood $V$ of $e$ such that $VX\subseteq XU$;

\vspace{3 mm}
$(iii)$   every bounded left uniformly continuous  function  on $G$ is right uniformly continuous. }
\vspace{5 mm}

A topological group $G$ is called {\it balanced} (or a  SIN-group)  if left and  right  uniformities  of $G$
coincide.
A group $G$  is called {\it functionally  balanced} if $G$  satisfies $(iii)$ of Theorem 3.4.
The study of functionally balanced groups  was initiated by  G. Itzkowitz [43].

The equivalence of $(ii)$  and $(iii)$  in Theorem 3.4 is a criterion for a topological group to be functionally balanced. In [44],  this criterion was used to show that each almost  metrizable  functionally  balanced  group is balanced.

\vspace{7 mm}

{\bf 3.3. Lattice topologies. } These topologies  on a complete lattice  $\mathcal{L}(G)$  of closed  subgroup are algebraically  defined by the lattice structure of $\mathcal{L}(G)$.

For example, a net $(A_{\alpha})_{\alpha\in\mathcal{I}}$  in   $\mathcal{L}(G)$
  {\it order-converges } to  $A\in  \mathcal{L}(G)$   if  there  exist two nets $(B_{\alpha})_{\alpha\in\mathcal{I}}$,
  $(C_{\alpha})_{\alpha\in\mathcal{I}}$ in $\mathcal{L}(G)$
   such that, for each
    $\alpha\in\mathcal{I}, $  $B_{\alpha}\subseteq A_{\alpha}\subseteq C_{\alpha}$
    and
    $\vee_{\alpha\in\mathcal{I}} B_{\alpha}= \wedge_{\alpha\in\mathcal{I}}C_{\alpha}=A $.
    By  [45],   for a compact group $G$,  every  net   in  $\mathcal{L}(G)$
       has  an order-convergent subset if and only if $\mathcal{L}(G)$
         endowed  with the  Shabauty  topology is a   topological  lattice, see  Theorem  1.10.

More on the lattices topologies on $\mathcal{L}(G)$  in the case of a compact $G$  can be find in  [46].

\vspace{7 mm}

{\bf 3.4.
Segment topologies. }
 Let $G$  be a topological  group,  $\mathcal{P}_{G}$  is the family of all subsets of $G$,
  $[G]^{<\omega}$    is the  family  of all finite subsets of  $G$.
Each  pair  $\mathcal{A}, \mathcal{B}$  of  subsets  of  $\mathcal{P}_{G}$  closed  under finite  unions define the segment  topology on  $\mathcal{L}(G)$  with a base  consisting of the segments.
$$[A, G\setminus B]= \{X\in\mathcal{L}(G): A\subseteq X\subseteq G\setminus B\}, \  A\in \mathcal{A}, \  B\in\mathcal{B}.$$

These topologies     are studied in [47]  in the following three cases:
$\mathcal{A}= \mathcal{B}= [G]^{<\omega}$;    $\mathcal{A}= \mathcal{P}_{G}$  and
$\mathcal{B}= [G]^{<\omega}$; $\mathcal{A}=[G]^{<\omega}, $   $\mathcal{B}= \mathcal{P}_{G}$

\vspace{7 mm}

{\bf 3.5.  $(\Sigma , \Theta )$-topologies.}
 This general  construction  for  topologizations  of the  set  $\mathcal{L}(G)$  of closed subgroups  of a  topological  group $G$  from [48]   produces Chabauty,  Vietoris, Bourbaki topologies and a plenty of  other topologies

We assume  that,  for each $H\in\mathcal{L}(G)$,  $\Sigma (H)$
  is some family  of open  subsets of $G$,
  $\Sigma=\cup _{H\in\mathcal{L}(G)}  \ \Sigma(H)$  and the following conditions are  satisfied

\begin{itemize}

\item{} if  $U, \mathcal{V}\in\Sigma(H)$ then  $U\cap\mathcal{V}$ contains some  $W\in\Sigma(H)$;

\item{} for every $U\in\Sigma(H)$, there exists $\mathcal{V}\in\Sigma(H)$ such that  $U\in\Sigma(K)$ for each
$K\in\mathcal{L}(G)$,  $K\subseteq\mathcal{V}$;

\item{} $\bigcap_{U\in\Sigma(H)} \overline{U}= H$  for each  $H\in\mathcal{L}(G).$

\end{itemize}

Then the family    $\{X\in \mathcal{L}(G): X\subseteq U\}$, $U\in \Sigma $   is a
base for the $\Sigma$-{\it topology} on  $\mathcal{L}(G)$.

Let $\tau$  denotes  the topology of $G $,  $\mathcal{P}_{\tau}$  is the family of all subsets of $\tau$.
We  assume that, for each  $H\in \mathcal{L}(G)$,   $\Theta(H)$  is some subset of  $\mathcal{P}_{\tau}$  such  that the  following conditions are satisfied

\begin{itemize}

\item{} for every $\alpha,\beta\in\Theta(H)$, there is  $\gamma\in\Theta(H)$  such that
$\alpha<\gamma$, $\beta<\gamma$ ($\alpha<\beta$ means that, for every $U\in\alpha$,  there  exists  $V\in \beta$   such that  $V\subseteq  U$);

\item{} for every $\alpha\in\Theta(H)$, there exists   $\beta\in\Theta(H)$  such that if
$K\in\mathcal{L}(G)$  and $K\cap V\neq\emptyset$ for each  $V\in \beta$, then $\alpha<\gamma$
for some $\gamma\in\Theta(K)$;

\item{} for each   $H\in\mathcal{L}(G)$ and every  neighbourhood   $V$ of $x$,  there exists $\alpha\in\Theta(H)$ such that $x\in U$,  $U\subseteq V$ for some  $U\in\alpha$.

\end{itemize}

Then the family    $\{X\in \mathcal{L}(G):X\cap  U\neq\emptyset$   for each   $U\in \alpha\}$,  where $\alpha\in\Theta(H)$, $H\in\mathcal{L}(G)$,  is a base for the $\Theta$-{\it topology}
  on $\mathcal{L}(G)$.

The upper bound of $\Sigma$- and  $\Theta$-topologies  is called the $(\Sigma, \Theta)$-{\it topology}.

A net  $(H_{\alpha})_{\alpha\in\mathcal{I}}$ converges in $(\Sigma, \Theta)$-topology
  to $H\in\mathcal{L}(G)$ if and only if

 \begin{itemize}

\item{} for any  $U\in\Sigma(H)$,  there exists $\beta\in\mathcal{I}$ such that $H_{\alpha}\subseteq U$ for each $\alpha>\beta$;

\item{}  for any  $\alpha\in\Theta(H)$,  there exists $\gamma\in\mathcal{I}$ such that $H_{\alpha}\cap \mathcal{V}\neq\emptyset$ for each $\alpha>\gamma$.

\end{itemize}

In [48],  one can find  characterizations  of $G$  with  compact and discrete
$\mathcal{L}(G)$  in some concrete  $(\Sigma, \Theta)$-topologies.

\vspace{7 mm}

{\bf 3.6.  Hyperballeans of groups.}
Let $G$  be a discrete  group.  The set
$\{Fg: g\in G,  F\in [G]^{<\omega}\}$
 is a family  of balls in the finitary coarse  structure  on $G$.
For coarse structures and balleans  see [49] and  [50].
The finitary  coarse  structure on $G$  induces the coarse  structure on   $\mathcal{L}(G)$ in which
$\{X\in\mathcal{L}(G): X\subseteq FA, \  A \in FX\}$,  $F\in [G]^{<\omega}$
  is the family  of balls centered  at  $A\in \mathcal{L}(G)$.  The set    $\mathcal{L}(G)$ endowed  with structure is called a hyperballean of  $G$.
Hyperballeans  of groups carefully  studied  in [51] can be considered  as asymptotic counterparts of Bourbaki  uniformities.


\vskip 5pt

CONTACT INFORMATION

I.~Protasov: \\
Faculty of Computer Science and Cybernetics  \\
        Kyiv University  \\
         Academic Glushkov pr. 4d  \\
         03680 Kyiv, Ukraine \\ i.v.protasov@gmail.com

\end{document}